\documentclass[11pt,reqno,a4paper]{amsart}

\usepackage[latin1]{inputenc}

\usepackage{etoolbox}
\usepackage{amsmath}
\usepackage{enumerate}
\usepackage{amssymb}
\usepackage{amscd}
\usepackage{amsthm}
\usepackage{amsfonts}
\usepackage{enumitem}

\patchcmd{\subsection}{-.5em}{.5em}{}{}
\patchcmd{\subsubsection}{-.5em}{.5em}{}{}

\usepackage{enumitem}

\usepackage[T1]{fontenc}

\usepackage[scaled]{helvet} 
\usepackage{courier} 
\usepackage{eulervm}
\normalfont

\makeatother
\usepackage{hyperref}

\numberwithin{equation}{section}

\newcommand{\cB}{\mathcal{B}}

\newcommand{\cH}{\mathcal{H}}

\newcommand{\cV}{\mathcal{V}}
\newcommand{\cW}{\mathcal{W}}

\newcommand{\bN}{\mathbb{N}}

\newcommand{\bZ}{\mathbb{Z}}

\newcommand{\ra}{\rightarrow}

\newcommand{\qand}{\quad \textrm{and} \quad}

\newcommand{\eps}{\varepsilon}

\DeclareMathOperator{\Hom}{Hom}

\DeclareMathOperator{\supp}{supp}
\DeclareMathOperator{\trace}{tr}

\DeclareMathOperator{\id}{id}
\DeclareMathOperator{\HS}{HS}

\DeclareMathOperator{\real}{Re}

\newtheorem{theorem}{Theorem}[section]

\newtheorem{proposition}[theorem]{Proposition}

\theoremstyle{definition}

\newtheorem{remark}[theorem]{Remark}

\renewcommand\labelenumi{(\roman{enumi})}
\renewcommand\theenumi\labelenumi

\usepackage{changepage}

\begin{document}

\title[Bohr sets in triple products]{Bohr sets in triple products of large sets in amenable groups}

\author{Michael Bj\"orklund}
\address{Department of Mathematics, Chalmers, Gothenburg, Sweden}
\email{micbjo@chalmers.se}

\author{John T. Griesmer}
\address{Department of Applied Mathematics and Statistics, Colorado School of Mines, Golden, Colorado}
\email{jtgriesmer@gmail.com}



\begin{abstract}
We answer a question of Hegyv\'ari and Ruzsa concerning effective estimates of the Bohr-regularity of
certain triple sums of sets with positive upper Banach densities in the integers. Our
proof also works for any discrete amenable group, and it does not require all addends in the triple products
we consider to have positive (left) upper Banach densities; one of the addends is allowed to only have positive
upper asymptotic density with respect to a (possibly very sparse) \emph{ergodic} sequence.
\end{abstract}

\date{}
\maketitle

\section{Main results}

Let $\Gamma$ be a discrete group. If $(\cW,\sigma)$ is a finite-dimensional unitary $\Gamma$-representation, we
define the (possibly trivial) semi-norm $|\cdot|_\sigma$ on $\Gamma$ by
\[
|\gamma|_\sigma = \sup\big\{ \|\sigma(\gamma)w-w\|_\cW \, : \, \|w\|_{\cW} = 1\big\}, \quad \textrm{for $\gamma \in \Gamma$},
\]
and if $S$ is a finite set of finite-dimensional \emph{irreducible} unitary $\Gamma$-representations and $\eps > 0$,
we define the \textbf{Bohr set} $U_{S,\eps}$ by
\[
U_{S,\eps} = \big\{ \gamma \in \Gamma \, : \, |\gamma|_\sigma < \eps, \enskip \textrm{for all $\sigma \in S$} \big\}.
\]
We denote by $\widehat{\Gamma}_{\textrm{FD}}$ the set of (unitary equivalence classes) of finite-dimensional irreducible
unitary $\Gamma$-representations. When $\Gamma$ is abelian, $\widehat{\Gamma}_{\textrm{FD}} \cong \Hom(\Gamma,S^1)$, where $S^{1}$ is the multiplicative group of complex numbers of modulus $1$; thus
\[
U_{S,\eps} = \big\{ \gamma \in \Gamma \, : \, |\xi(\gamma) -1| < \eps, \enskip \textrm{for all $\xi \in S$} \big\},
\]
where $S$ is a finite subset of $\Hom(\Gamma,S^1)$. If $\Gamma$ lacks non-trivial finite-dimensional unitary representations, for instance, if $\Gamma$ is a finitely generated infinite simple group, then the only Bohr set is $\Gamma$ itself. \\

Let $B$ be a subset of $\Gamma$. Given a sequence $(T_n)$ of finite sets in $\Gamma$,
we define the \textbf{upper asymptotic density} of $B$ along $(T_n)$ by
\begin{equation}
\label{eq_asymdens}
\overline{d}_{(T_n)}(B) = \varlimsup_n \, \frac{|B \cap T_n|}{|T_n|},
\end{equation}
and if $\Gamma$ is amenable, we define the \textbf{upper Banach density} of $B$ by
\begin{equation}
\label{eq_UBD}
d^*(B) = \sup\big\{ \overline{d}_{(F_n)}(B) \, : \, \textrm{$(F_n)$ is a left F\o lner sequence} \big\}.
\end{equation}
A sequence $(T_n)$ of finite subsets of $\Gamma$ is called \textbf{ergodic} (or equidistributed) if for every unitary
$\Gamma$-representation $(\cH,\pi)$ and for all $u \in \cH$,
\[
\frac{1}{|T_n|} \sum_{\gamma \in T_n} \pi(\gamma)u \ra P_\pi u,
\]
in the weak topology on $\cH$, where $P_\pi$ denotes the orthogonal projection onto the subspace $\cH^\pi$ consisting of $\pi(\Gamma)$-fixed points in $\cH$.

By the Weak Ergodic Theorem, every left F{\o}lner sequence in a discrete amenable group is ergodic. We stress that the
converse is far from true, ergodic sequences can be very sparse; for instance, the sequence $T_n = \big\{ \lfloor k^{5/2} \rfloor \, : \, k = 1,\ldots,n \big\}$ is ergodic in $\Gamma = (\bZ,+)$ (see for instance \cite{BKQW}), while $\bigcup_{n\in \mathbb N} T_{n}$ has zero upper Banach density. \\

It is a fundamental problem in arithmetic combinatorics to determine which products of "large" subsets of $\Gamma$ contain Bohr sets. The story here starts with the influential work of
Bogolio\`uboff \cite{Bog}, who proved that for any subset $B \subset \bZ$ with positive upper Banach density, the
four-fold difference set $B-B+B-B$ contains a Bohr set $U_{S,\eps}$, whose parameters $|S|$ and $\eps$ can be estimated
in terms of $d^*(B)$ alone. A non-effective generalization was later provided by F\o lner \cite{Fol}, who showed that if $\Gamma$
is a discrete amenable group and $B \subset \Gamma$ has positive upper Banach density, then the difference set $BB^{-1}$ contains the intersection of a Bohr set $U_{S,\eps}$ and a subset $T \subset \Gamma$, whose complement has zero upper Banach density.
The question of whether the set $T$ could be removed from the conclusion attracted quite a lot of attention until K\v{r}\'{\i}\v{z}
\cite{Kriz} finally answered it in the negative.

Before K\v{r}\'{\i}\v{z}'s result, Ellis and Keynes proved a result between Bogolio\`uboff's and F\o lner's, namely that if $\Gamma$ is a
discrete amenable group and $B \subset \Gamma$ is right syndetic (meaning that there is a finite set $F \subset \Gamma$
such that $FB = \Gamma$), then there exists a Bohr set $U_{S,\eps}$ such that the triple product
$BBB^{-1}$ contains a set of the form $b_o U_{S,\eps}$ for some $b_o \in B$. However, their method does not bound the parameters $|S|$ and $\eps$ in terms of the size of $B$. More recently,
Hegyv\'ari and Ruzsa \cite{HR} established the same result for $\Gamma = (\bZ,+)$, but assuming only that $B$ has positive
upper Banach density.   Section 3 of \cite{HR} explains why the method therein provides no bounds on $|S|$ and $\eps$ in terms $d^*(B)$, and asks if this can be done. The aim of  this paper is to provide such bounds, in a more general setting. \\

Our first main result can be stated as follows.

\begin{theorem}
\label{thm_maincomb}
Let $\Gamma$ be a discrete amenable group and let $(T_n)$ be an ergodic sequence in $\Gamma$.
Suppose that $A$ and $B$ are subsets of $\Gamma$ with
\[
\overline{d}_{(T_n)} =\alpha > 0 \qand d^*(B) = \beta > 0.
\]
Then there exist $a_o \in A$ and a set $S \subset \widehat{\Gamma}_{\textsc{FD}} \setminus \{[\id]\}$ with $|S| \leq \frac{17(1-\alpha^2)}{d_\Gamma \alpha^2 \beta^2}$ such that
\[
ABB^{-1} \supset a_o \, U_{S,\frac{\beta}{2}},
\]
where $d_\Gamma$ denotes the smallest dimension of a non-trivial unitary $\Gamma$-representation. If $d_\Gamma = \infty$,
that is to say, if $\Gamma$ does not admit any non-trivial finite-dimensional unitary representations, or if $\alpha = 1$, then
we have $ABB^{-1} = \Gamma$.
\end{theorem}

\begin{remark}
Theorem \ref{thm_maincomb} applies to the setting when $A$ has positive upper Banach density $\alpha$, and $(T_n)$ is
a left F\o lner sequence (and thus an ergodic sequence) in $\Gamma$ such that $\overline{d}_{(T_n)}(A) = d^*(A) = \alpha$; the  special case where $\Gamma = (\bZ,+)$ and $A = B$ provides the bounds requested in \cite{HR}.
\end{remark}

Let us now discuss a "dynamical" strengthening of Theorem \ref{thm_maincomb}. First note that the triple product set
$ABB^{-1}$ consists of exactly those $\gamma$ in $\Gamma$ such that $AB \cap \gamma B \neq \emptyset$. In what
follows we estimate the content of these intersections. We shall show that for any number
$0 \leq r < d^*(B)$, the level set
\[
\big\{ \gamma \in \Gamma \, : \, d^*(AB \cap \gamma B) > d^*(B) r \big\}
\]
contains a translate of a Bohr set $U_{S,\eps}$ whose parameters $|S|$ and $\eps$ can be estimated in terms of $d^*_{(T_n)}(A)$, $d^*(B)$ and $r$ alone. The arguments needed to deduce such a statement from our next result are by now standard, and will be omitted here; the reader could for instance consult Section 2 in \cite{Bj}.  Note that Theorem \ref{thm_mainerg} is valid for any discrete group $\Gamma$, amenable or not.

\begin{theorem}
\label{thm_mainerg}
Let $\Gamma$ be a discrete group, let $(T_n)$ be an ergodic sequence in $\Gamma$ and let $(Y,\nu)$ be a probability
measure preserving $\Gamma$-space. Suppose that $A$ is a subset of $\Gamma$ and $B$ is a Borel subset of $Y$ with
\[
\overline{d}_{(T_n)}(A) = \alpha > 0 \qand \nu(B) = \beta > 0.
\]
Then, for every $0 \leq r < \beta$, there exist $a \in A$ and $S \subset \widehat{\Gamma}_{\textsc{FD}} \setminus \{[id]\}$ with
$|S| \leq \frac{17(1-\alpha^2)}{d_\Gamma \alpha^2 (\beta-r)^2}$ such that
\[
\big\{ \gamma \in \Gamma \, : \, \nu(AB \cap \gamma B) > \beta r \big\}
\supset
a \, U_{S,\frac{\beta-r}{2}},
\]
where $d_\Gamma$ denotes the smallest dimension of a non-trivial unitary $\Gamma$-representation.
If $d_\Gamma = \infty$, that is to say, if $\Gamma$ does not admit any non-trivial finite-dimensional unitary representations,
or if $\alpha = 1$, then $\nu(AB) = 1$.
\end{theorem}

\begin{remark}
Theorem \ref{thm_mainerg} is in stark contrast with Corollary 2.6 in \cite{Gr1} by the second author, which roughly
asserts that for $\Gamma = (\bZ,+)$, one can construct an ergodic probability measure-preserving $\Gamma$-space $(Y,\nu)$ and
a Borel set $B \subset Y$ with $\nu$-measure arbitrarily close to $\frac{1}{2}$, such that \emph{no} level set of the form
$\{ \gamma \, : \, \nu(B \cap \gamma B) > t \big\}$ for $t > 0$ contains a translate of a Bohr set. In \cite{Gr2}, the second author provides an even
stronger counterexample for $\Gamma = (\bigoplus_{\bZ} \bZ/2\bZ,+)$.
\end{remark}

We will deduce Theorem \ref{thm_mainerg} from a general result about unitary $\Gamma$-representations. Before we state this result, we define a slight generalization of the condition that $A$ has positive density along some ergodic sequence,
introducing some notation and conventions on the way. Let $(T_n)$ be an ergodic sequence in a discrete group $\Gamma$ and
let $A$ be a subset of $\Gamma$ with $\overline{d}_{(T_n)}(A) = \alpha > 0$. We consider the sequence $(\lambda_n)$ of means (positive and unital functionals on $\ell^\infty(\Gamma)$) defined by
\begin{equation}
\label{eq_fromTntoFS}
\lambda_n(f) = \frac{1}{|T_n|} \sum_{\gamma \in T_n} f(\gamma), \quad \textrm{for $f \in \ell^\infty(\Gamma)$}.
\end{equation}
and pick a weak*-cluster point $\lambda$ such that $\alpha = \lambda(\chi_A)$. Recall that the \textbf{Fourier-Stiltjes algebra}
$\cB(\Gamma)$ is the $*$-algebra of bounded functions on $\Gamma$ spanned by all matrix coefficients of unitary $\Gamma$-representations. It is not hard to prove
(see for instance \cite{God}) that there exists a \emph{unique} left $\Gamma$-invariant mean $\eta$ on $\cB(\Gamma)$. Since
$(T_n)$ is ergodic, we see that $\lambda|_{\cB(\Gamma)} = \eta$. Any mean on $\ell^\infty(\Gamma)$ with this property
will be called a \textbf{Fourier-Stiltjes} mean.

In what follows, we shall adopt two convenient abuses of notation concerning means on $\ell^\infty(\Gamma)$. Firstly, if $A \subset \Gamma$ and $\lambda$ is a mean, we write $\lambda(A) = \lambda(\chi_A)$, where $\chi_A$ denotes the indicator function on $A$. Secondly, for $f \in \ell^\infty(\Gamma)$, we write
\[
\int_A f \, d\lambda = \lambda(\chi_A f),
\]
although the expression on the left hand side is not an integral in the Lebesgue sense. We can now formulate the main result of
this paper.

\begin{theorem}
\label{thm_mainuni}
Let $\Gamma$ be a discrete group, $(\cH,\pi)$ a unitary $\Gamma$-representation and let $\lambda$ be a Fourier-Stiltjes mean
on $\Gamma$. Suppose that $A$ is a subset of $\Gamma$ with $\alpha = \lambda(A) > 0$. Then, for every $\beta > 0$ and for
every $0 \leq r < \beta$, there exists a subset $S \subset \widehat{\Gamma}_{\textrm{FD}}$ with $|S| \leq \frac{17(1-\alpha^2)}{d_\Gamma \alpha^2(\beta-r)^2}$ such that for any unit vector $u \in \cH$ with $\|P_\pi u\|_{\cH} \geq \sqrt{\beta}$, we have
\begin{equation}
\label{eq_generalincl}
\big\{ \gamma \in \Gamma \, : \, \frac{1}{\lambda(A)} \int_{A} \real \big\langle \pi(a)u, \pi(\gamma)u \big\rangle_{\cH} \, d\lambda(a) > r \big\}
\supset a_o U_{S, \frac{\beta - r}{2}},
\end{equation}
for some $a_o \in A$, where $P_\pi$ denotes the orthogonal projection onto $\cH^{\pi}$.
\end{theorem}

\subsection{Organization of the paper}

In Section \ref{sec:mainunitomainerg} we show how Theorem \ref{thm_mainerg} can be deduced from Theorem \ref{thm_mainuni}.
In Section \ref{sec:outline} we break down the proof of Theorem \ref{thm_mainuni} into two main propositions, which will be proved
in Section \ref{sec:main1} and Section \ref{sec:main2} respectively.

\section{Proof of Theorem \ref{thm_mainerg} assuming Theorem \ref{thm_mainuni}}
\label{sec:mainunitomainerg}

Let $\Gamma$ be a discrete group and let $(Y,\nu)$ be a probability-measure preserving $\Gamma$-space. The regular
(unitary) $\Gamma$-representation $(L^2(\nu),\pi)$ is defined by
\[
(\pi(\gamma)f)(y) = f(\gamma^{-1}y), \quad \textrm{for $f \in L^2(\nu)$}.
\]
Let $B \subset Y$ be a Borel set with $\beta = \nu(B) > 0$. Then $u_B = \frac{\chi_B}{\nu(B)^{1/2}}$ is a unit vector in $L^2(\nu)$,
and since $1 \in L^2(\nu)^\pi$ is a unit vector as well, we have
\[
\|P_\pi u_B \|_{L^2(\nu)} \geq  \big\langle u_B, 1 \big\rangle_{L^2(\nu)} = \sqrt{\beta}.
\]
Let $(T_n)$ be an ergodic sequence in $\Gamma$ and suppose that $A$ is a subset of $\Gamma$ with
$\overline{d}_{(T_n)}(A) = \alpha > 0$. Then the construction in \eqref{eq_fromTntoFS} produces a Fourier-Stiltjes mean
$\lambda$ on $\Gamma$ such that $\lambda(A) = \alpha$. Note that
\begin{eqnarray*}
\frac{\nu(AB \cap \gamma B)}{\nu(B)}
&\geq &
\frac{1}{\lambda(A)} \int_A \frac{\nu(aB \cap \gamma B)}{\nu(B)} \, d\lambda(a) \\
&=&
\frac{1}{\lambda(A)} \int_A \big\langle \pi(a) u_B, \pi(\gamma)u_B \big\rangle_{L^2(\nu)} \, d\lambda(a),
\end{eqnarray*}
and thus for all $r \geq 0$,
\[
\big\{ \gamma \in \Gamma \, : \, \nu(AB \cap \gamma B) > \beta r \big\}
\supset
\Big\{ \gamma \in \Gamma \, : \, \frac{1}{\lambda(A)} \int_A \big\langle \pi(a) u_B, \pi(\gamma)u_B \big\rangle_{L^2(\nu)} \, d\lambda(a) > r \Big\}.
\]
Theorem \ref{thm_mainuni} tells us that for every $0 \leq r < \beta$, there exists $a_o \in A$ such that the set on the right hand side
contains a set of the form $a_o U_{S,\frac{\beta-r}{2}}$, for some $S$ with $|S| \leq \frac{17(1-\alpha^2)}{d_\Gamma \alpha^2 (\beta-r)^2}$.

\section{An outline of the proof of Theorem \ref{thm_mainuni}}
\label{sec:outline}

The aim of this section is to break up the proof of Theorem \ref{thm_mainuni} into two main propositions which will be proved
in Section \ref{sec:main1} and Section \ref{sec:main2} respectively. Once we have stated these propositions and provided some background, we show how Theorem \ref{thm_mainuni} can be deduced from these results.\\

Let $\Gamma$ be a discrete group. A \textbf{compactification} of $\Gamma$ is a pair $(K,\tau)$, where $K$ is a compact
Hausdorff group and $\tau : \Gamma \ra K$ is a (not necessarily injective) homomorphism. We denote by $m_K$ the (unique)
Haar probability measure on $K$, and by $\widehat{K}$ the set of equivalence classes of strongly continuous \emph{irreducible}
unitary $K$-representations. \\

If $(\cW,\sigma)$ is a strongly
continuous irreducible unitary $K$-representation, we define the continuous semi-norm $|\cdot|_\sigma$ on $K$ by
\begin{equation}
\label{def_normsigma}
|k|_\sigma = \sup\{ \|\sigma(k)w-w\|_{\cW} \, : \, \|w\|_{\cW} = 1 \big\},
\end{equation}
and for a (possibly infinite) subset $S$ of inequivalent strongly continuous irreducible unitary $K$-representations, we
define
\begin{equation}
\label{def_normsigmaS}
|k|_S = \sup\{ |k|_\sigma \, : \, \sigma \in S \big\}, \quad \textrm{for $k \in K$}.
\end{equation}  
The inequality
\begin{equation}\label{ineq_subadditive}
  |k_{1}k_{2}|_{S}\leq |k_{1}|_{S}+|k_{2}|_{S} \quad \textrm{for all $k_{1}, k_{2}\in K$}
\end{equation}
follows immediately from the definition and subadditivity of the Hilbert space norm.

Note that we can only ensure that $|\cdot|_S$ is continuous if $S$ is \emph{finite}. If $\mu$ is a Borel probability measure
on $K$ and $(\cV,\theta)$ is a strongly continuous (not necessarily irreducible) unitary $K$-representation, we define a linear
operator $\theta(\mu) : \cV \ra \cV$ by
\[
\theta(\mu)v = \int_K \theta(k) v \, d\mu(k), \quad \textrm{for $v \in \cV$}.
\]
Among all compactifications of $\Gamma$, there is a universal one, known as the \textbf{Bohr compactification}, which
here will be denote by $(b\Gamma,\iota)$. One way to construct it is by taking the closure of $(\tau_{i}(\Gamma))_{i\in E}$ in the direct product $\prod_{i \in E} K_{i}$, where $\{(K_{i},\tau_{i}):i\in E\}$ is the set of (representatives of) \emph{all} compactifications of $\Gamma$ indexed by a set $E$; we stress that such a product is not necessarily metrizable (see Section 4.7 in \cite{Foll} for a more detailed discussion about the Bohr compactification for abelian groups). Universality of $(b\Gamma,\iota)$ now amounts to the following
property: Whenever $K$ is a compact Hausdorff group and $\tau : \Gamma \ra K$ is a homomorphism, then there exists
a unique \emph{continuous} homomorphism $\overline{\tau} : b\Gamma \ra K$ such that $\tau = \overline{\tau} \circ \iota$. \\

With these definitions we can break Theorem \ref{thm_mainuni} into two technical propositions.

\begin{proposition}
\label{prop_main1}
Let $\Gamma$ be a discrete group, $(\cH,\pi)$ a unitary $\Gamma$-representation and let $\lambda$ be a Fourier-Stiltjes mean
on $\Gamma$. Suppose that $A$ is a subset of $\Gamma$ with $\alpha = \lambda(A) > 0$. Then there exist
\begin{enumerate}
\item a regular Borel probability measure $\mu_A$ on $b\Gamma$, absolutely continuous with respect to $m_{b\Gamma}$,
such that
\[
\supp \mu_A \subset \overline{\iota(A)} \qand \big\|\frac{d\mu_A}{dm_{b\Gamma}}\big\|_{L^\infty(b\Gamma)} \leq \frac{1}{\alpha},
\]
\item a strongly continuous unitary $b\Gamma$-representation $(\cV,\theta)$ and a linear $\Gamma$-equivariant map
\[
Q : (\cH,\pi) \ra (\cV,\theta \circ \iota)
\]
with $\|Q\|_{\textrm{op}} \leq 1$, whose restriction $Q : \cH^\pi \ra \cV^{\theta \circ \iota}$ is an isometry,
\end{enumerate}
such that
\[
\real \big\langle \theta(\mu_A)Qu,\theta(\iota(\gamma))Qu \big\rangle_{\cV}
=
\frac{1}{\lambda(A)} \int_{A} \real \big\langle \pi(a)u, \pi(\gamma)u \big\rangle_{\cH} \, d\lambda(a),
\]
for all $u \in \cH$ and $\gamma \in \Gamma$.
\end{proposition}

\begin{proposition}
\label{prop_main2}
Let $K$ be a compact Hausdorff group, $\mu$ a regular Borel probability measure on $K$ and let $(\cV,\theta)$ be a strongly
continuous unitary $K$-representation. Fix $\delta > 0$ and set
\[
S_\delta = \big\{ [\sigma] \in \widehat{K} \setminus \{[\id]\} \, : \, \|\sigma(\mu)\|_{\textrm{op}} > \delta \big\}.
\]
Then, for every $v \in \cV$, there exists $k_o \in \supp \mu$ such that
\begin{equation}
\label{eq_cptincl}
\real \big\langle \theta(\mu)v,\theta(k_o k)v \big\rangle_{\cV} \geq \|P_\theta v\|^2_{\cV} - (2\delta + |k|_{S_\delta})\|v\|^2_{\cV},
\end{equation}
for all $k \in K$. Furthermore, if $\mu \ll m_K$ and $\|\frac{d\mu}{dm_K}\|_{L^2(K)} \leq \xi$, then
\begin{equation}
\label{eq_uppbndSd}
|S_\delta| \leq \frac{\xi^2-1}{d_K \delta^2}.
\end{equation}
\end{proposition}

\subsection{Proof of Theorem \ref{thm_mainuni} assuming Proposition \ref{prop_main1} and Proposition \ref{prop_main2}}

Let $\Gamma$ be a discrete group, $(\cH,\pi)$ a unitary $\Gamma$-representation and let $\lambda$ be a Fourier-Stiltjes mean
on $\Gamma$. Suppose that $A$ is a subset of $\Gamma$ with $\alpha = \lambda(A) > 0$. Fix $\beta > 0$ and suppose that
$u$ is a unit vector in $\cH$ with $\|P_\pi u\| \geq \sqrt{\beta}$. Let us also fix $0 \leq r < \beta$. Let $\mu_A$, $(\cV,\theta)$ and
$Q$ be as in Proposition \ref{prop_main1}, so that
\begin{equation}
\label{eq_passage}
\real \big\langle \theta(\mu_A)v,\theta(\iota(\gamma))v \big\rangle_{\cV}
=
\frac{1}{\lambda(A)} \int_{A} \real \big\langle \pi(a)u, \pi(\gamma)u \big\rangle_{\cH} \, d\lambda(a),
\end{equation}
with $v = Qu$, for all $\gamma \in \Gamma$. We stress that $\mu_A$ does not depend on the choice of $u$. Fix $\delta > 0$ and
let $S_\delta$ be as in Proposition \ref{prop_main2}, applied to $K = b\Gamma$, so that for some $k_o \in \supp(\mu_A)$,
\begin{equation}
\label{eq_lowbndproof}
\real \big\langle \theta(\mu)v,\theta(k_o k)v \big\rangle_{\cV} \geq \|P_\theta v\|^2_{\cV} - (2\delta + |k|_{S_\delta})\|v\|^2_{\cV},
\end{equation}
for all $k \in b\Gamma$. Since $\big\|\frac{d\mu_A}{dm_{b\Gamma}}\big\|_{L^2(b\Gamma)} \leq \frac{1}{\alpha}$, it follows from
the second part of Proposition \eqref{prop_main2} that
\begin{equation}
\label{eq_Sdelta}
|S_\delta| \leq \frac{1-\alpha^2}{d_\Gamma \alpha^2 \delta^2}.
\end{equation}
Furthermore, since $\|Q\|_{\textrm{op}} \leq 1$ and the restriction of $Q$ to $\cH^\pi$ is an isometry, we see that
\[
\|v\| = \|Qu\| \leq 1 \qand \|P_\theta v\|_{\cV} = \|P_\theta Q u\|_{\cV} = \|Q P_\pi u\|_{\cV} = \|P_\pi u\|_{\cH} \geq \sqrt{\beta}.
\]
Plugging this into \eqref{eq_lowbndproof}, we get
\begin{equation}
\label{eq_almostfinallowbnd}
\real \big\langle \theta(\mu)v,\theta(k_o k)v \big\rangle_{\cV} \geq \beta - (2\delta + |k|_{S_\delta}), \quad \textrm{for all $k \in b\Gamma$}.
\end{equation}
Choose $T > 0$ such that
\begin{equation}
\label{eq_T}
\big(1 - \frac{1}{2T}\big)\big(1 + \frac{1}{2T}\big)^2 > 1 \qand \big(4 + \frac{2}{T}\big)^2 \leq 17
\end{equation}
and set
\[
\delta = \frac{\beta - r}{4 + \frac{2}{T}} > 0.
\]
By \eqref{eq_Sdelta}, we get the bound
\[
|S_\delta| \leq \frac{\big(4 + \frac{2}{T}\big)^2 (1-\alpha^2)}{d_\Gamma \alpha^2 (\beta-r)^2} \leq
\frac{17(1-\alpha^2)}{d_\Gamma \alpha^2 (\beta-r)^2}.
\]
Since $k_o \in \supp \mu_A \subset \overline{\iota(A)}$ and $S_\delta$ is finite, we can find $a_o \in A$ such that $|k_o^{-1}\iota(a_o)|_{S_\delta} \leq \frac{\beta-r}{4T^2}$, and thus, for all $k \in b\Gamma$ with $|k|_{S_\delta} \leq \frac{\beta-r}{2}$, we have by \eqref{eq_almostfinallowbnd} and \eqref{ineq_subadditive},
\begin{eqnarray*}
\real \big\langle \theta(\mu)v,\theta(\iota(a_o)k)v \big\rangle_{\cV}
&=&
\real \big\langle \theta(\mu)v,\theta(k_o k_o^{-1}\iota(a_o) k)v \big\rangle_{\cV} \\
&\geq &
\beta - (2\delta + |k|_{S_\delta} + |k_o^{-1}\iota(a_o)|_{S_\delta}) \\
&\geq &
\frac{\beta}{2} + \frac{r}{2} - \Big(2\delta + \frac{\beta-r}{4T^2}\Big) > r,
\end{eqnarray*}
where the last inequality holds since
\[
1 > \frac{1}{1 + \frac{1}{2T}} + \frac{1}{(2T)^2}
\]
by the first condition in \eqref{eq_T}. Combining this with \eqref{eq_passage}, we conclude that
\[
\frac{1}{\lambda(A)} \int_{A} \real \big\langle \pi(a)u, \pi(a_o \gamma)u \big\rangle_{\cH} \, d\lambda(a) > r
\]
for all $\gamma \in \Gamma$ such that $|\gamma|_{S_\delta} \leq \frac{\beta-r}{2}$ (where we now consider $|\cdot|_{S_\delta}$
as a semi-norm on $\Gamma$ via the pull-back under $\iota$).

\section{Proof of Proposition \ref{prop_main1}}
\label{sec:main1}

Let $\Gamma$ be a discrete group, $(\cH,\pi)$ a unitary $\Gamma$-representation and let $\lambda$ be a Fourier-Stiltjes mean on
$\Gamma$. Let $\cH_o$ denote the
closure of the linear span of all finite-dimensional sub-$\Gamma$-representations of $(\cH,\pi)$, and set $\cH_1 = \cH_o^{\perp}$
so that $\cH \cong \cH_o \widehat{\oplus} \cH_1$. It is well-known (see for instance Chapter 2 in \cite{KeLi}) that
\[
\int_\Gamma \big| \big\langle \pi(a)u_1,u_1' \big\rangle_{\cH} \big|^2 \, d\eta(a) = 0, \quad \textrm{for all $u_1, u_1' \in \cH_1$},
\]
where $\eta$ denotes the unique $\Gamma$-invariant mean on $\cB(\Gamma)$, and thus for all $u_1 \in \cH_1$,
\[
\big| \int_A  \big\langle \pi(a)u_1,\pi(\gamma)u_1 \big\rangle_{\cH}  \, d\lambda(a) \big|
\leq
\Big( \int_\Gamma  \big| \big\langle \pi(a)u_1,\pi(\gamma)u_1 \big\rangle_{\cH} \big|^2  \, d\eta(a) \Big)^{1/2}
= 0,
\]
for all $\gamma \in \Gamma$ and $A \subset \Gamma$, since $\lambda$ restricts to $\eta$ on $\cB(\Gamma)$. In particular,
for every $A \subset \Gamma$,
\[
\int_A   \big\langle \pi(a)u,\pi(\gamma)u \big\rangle_{\cH} \, d\lambda(a)
=
\int_A  \big\langle \pi(a)P_o u,\pi(\gamma)P_o u \big\rangle_{\cH} \, d\lambda(a),
\]
for all $\gamma \in \Gamma$ and $u \in \cH$, where $P_o$ denotes the orthogonal projection of $\cH$ onto $\cH_o$. Since
$\cH_o$ is the closure of the linear span of the finite-dimensional unitary sub-representations of $\cH$, we note that the
restriction of $\pi$ to $\cH_o$ defines a homomorphism of $\Gamma$ into the direct product $K$ of the (compact) unitary groups
of the finite-dimensional sub-representations of $\cH$, which is again a compact Hausdorff group. By the universal property of the Bohr
compactification $(b\Gamma,\iota)$, this homomorphism extends to a continuous homomorphism $\theta$ from $b\Gamma$ into $K$. We conclude that
\[
\pi \mid_{\cH_o} = \theta \circ \iota,
\]
and thus
\begin{equation}
\label{eq_frompitotheta}
\int_A   \big\langle \pi(a)u,\pi(\gamma)u \big\rangle_{\cH} \, d\lambda(a)
=
\int_A  \big\langle \theta(\iota(a))P_o u,\theta(\iota(\gamma))P_o u \big\rangle_{\cH_o} \, d\lambda(a).
\end{equation}
In what follows, let $(\cV,\theta) = (\cH_o,\pi_o)$, which we regard as a strongly continuous unitary $b\Gamma$-representation,
and we denote by $Q$ the projection $P_o : \cH \ra \cH_o$. Since $\cH^\pi \subset \cH_o$, we see that $Q$ restricts to
the identity (hence an isometry) on $\cH^\pi$.\\

Let us now fix a subset $A \subset \Gamma$ with $\alpha = \lambda(A) > 0$, and define a regular Borel probability measure
$\mu_A$ on $b\Gamma$ by
\[
\mu_A(f) = \frac{1}{\lambda(A)} \int_A f(\iota(\gamma)) \, d\lambda(\gamma), \quad \textrm{for $f \in C(b\Gamma)$}.
\]
Since $\lambda$ is a Fourier-Stiltjes mean, we claim that $\iota_*\lambda = m_{b\Gamma}$. Indeed, the sub-$*$-algebra $\iota^*(C(b\Gamma)) \subset \ell^\infty(\Gamma)$, also known as the space of almost periodic functions on $\Gamma$, is contained in the Fourier-Stiltjes
algebra of $\Gamma$ (see for instance Theorem 1.43(5) in \cite{Glas}), and thus $\iota_*\eta = m_{b\Gamma}$ by the uniqueness
of $\iota(\Gamma)$-invariant probability measures on $b\Gamma$. Since $\lambda$ restricts to $\eta$ on $\iota^*(C(b\Gamma))$,
the claim follows. In particular,
\[
|\mu_A(f)| \leq \frac{1}{\alpha} \int_\Gamma |(f \circ \iota)| \, d\lambda = \frac{1}{\alpha} \|f\|_{L^1(b\Gamma)},
\]
for all $f \in C(b\Gamma)$. This allows us by Hahn-Banach's Theorem to extend $\mu_A$ to a linear functional on $L^1(b\Gamma)$
with the same norm, so we conclude that there exists $\rho_A \in L^\infty(b\Gamma)$ with $\|\rho\|_\infty \leq \frac{1}{\alpha}$ such that
\[
\mu_A(f) = \int_{b\Gamma} f(k) \overline{\rho_A(k)} \, dm_{b\Gamma}(k), \quad \textrm{for all $f \in C(b\Gamma)$}.
\]
It readily follows that $\rho_A$ is essentially non-negative with $m_{b\Gamma}$-integral one. Hence, $\mu_A$ is absolutely continuous
with respect to $m_{b\Gamma}$ and $\rho_A = \frac{d\mu_A}{dm_{b\Gamma}}$. Note that $\|\rho_A\|_{L^2(K)} \leq \|\rho_A\|_{L^\infty(K)} \leq \frac{1}{\alpha}$. Finally, pick $k \in \supp \mu_A$ and an open neighborhood $V$ of $k$ in $b\Gamma$. Let $f$ be
a non-negative, but not identically zero, continuous function on $b\Gamma$ which is supported inside $V$. Then
\[
0 < \mu_A(f) = \frac{1}{\lambda(A)} \int_A f(\iota(a)) \, d\lambda(a),
\]
which shows that there exists at least one $a \in A$ such that $\iota(a) \in V$. By considering smaller and smaller open neighborhoods $V$ around $k$, we can conclude that $k \in \overline{\iota(A)}$. Since $k$ is arbitrary, $\supp \mu_A \subset \overline{\iota(A)}$. \\

If we apply the argument above to the continuous functions
\[
k \mapsto \big\langle \theta(k)Qu,\theta(\iota(\gamma))Qu \big\rangle_{\cV}, \quad \textrm{for $k \in b\Gamma$ and $\gamma \in \Gamma$},
\]
then we see that
\begin{eqnarray*}
\frac{1}{\lambda(A)} \int_A \big\langle \theta(\iota(a))Qu,\theta(\iota(\gamma))Qu \big\rangle_{\cV} \, d\lambda(a)
&=&
\int_K \big\langle \theta(k)Qu,\theta(\iota(\gamma))Qu \big\rangle_{\cV} \, d\mu_A(k) \\
&=&
\big\langle \theta(\mu_A)Qu,\theta(\iota(\gamma))Qu \big\rangle_{\cV},
\end{eqnarray*}
for all $\gamma \in \Gamma$, which in combination with \eqref{eq_frompitotheta} finishes the proof.

\section{Proof of Proposition \ref{prop_main2}}
\label{sec:main2}

Let $K$ be a compact Hausdorff group, $(\cV,\theta)$ a strongly continuous unitary $K$-representation and let $\mu$ be a regular
Borel probability measure on $K$.
By Theorem 15.1.3 in \cite{Dix} (or Section 5.3 in \cite{Foll}), there exist a (possibly uncountable, if $K$ is not second countable) subset $\Omega$ of (unitarily inequivalent) strongly continuous \emph{irreducible} (hence
finite-dimensional) unitary $K$-representations $(\cW_\sigma,\sigma)$, indexed by $\sigma \in \Omega$, and a
function $m : \Omega \ra \bN \cup \{\infty\}$ such that
\[
\cV \cong \bigoplus_{\sigma \in \Omega} \cW_\sigma^{\oplus m_\sigma}.
\]
In particular, for every $v \in \cV$, there exist, for every $\sigma \in \Omega$, a vector $w_\sigma \in \cW_\sigma^{\oplus m_\sigma}$
and mutually orthogonal vectors $w_{\sigma,j} \in \cW_\sigma$ for $j=1,\ldots,m_\sigma$, such that
\[
v = \sum_{\sigma} w_\sigma \qand \|v\|_{\cV}^2 = \sum_{\sigma} \|w_\sigma\|_{\cW_\sigma^{\oplus m_\sigma}}^2 \qand
\|w_\sigma\|_{\cW_\sigma^{\oplus m_\sigma}}^2
= \sum_{j=1}^{m_\sigma} \|w_{\sigma,j}\|^2_{\cW_\sigma},
\]
for every $\sigma \in \Omega$. Hence, for every $k \in K$,
\[
\big\langle \theta(\mu)v,\theta(k)v \big\rangle_{\cV}
=
\sum_{\sigma} \Big( \sum_{j=1}^{m_\sigma} \big\langle \sigma(\mu)w_{\sigma,j},\sigma(k)w_{\sigma,j} \big\rangle_{\cW_\sigma} \Big),
\]
and thus,
\begin{eqnarray}
\big|
\big\langle \theta(\mu)v,\theta(k)v \big\rangle_{\cV}
-
\big\langle \theta(\mu)v,\theta(l)v \big\rangle_{\cV}
\big|
&\leq &
\sum_{\sigma} \Big( \sum_{j=1}^{m_\sigma} \|\sigma(\mu)\|_{\textrm{op}} |k^{-1}l|_{\sigma} \|w_{\sigma,j}\|_{\cW_\sigma}^2 \Big) \nonumber \\
&=&
\sum_{\sigma} \|\sigma(\mu)\|_{\textrm{op}} |k^{-1}l|_{\sigma} \|w_{\sigma}\|_{\cW_\sigma^{\oplus m_\sigma}}^2,
\label{lowbnd}
\end{eqnarray}
for all $k,l \in K$, where $|\cdot|_\sigma$ is defined as in \eqref{def_normsigma}. Fix $\delta > 0$ and set
\[
S_\delta = \big\{ \sigma \in \Omega \setminus \{\id\} \, : \, \|\sigma(\mu)\|_{\textrm{op}} > \delta \big\},
\]
where $\id$ denotes the trivial $K$-representation. Then, \eqref{lowbnd}, combined with the fact that
\[
|k|_\sigma \leq 2 \qand \|\sigma(\mu)\|_{\textrm{op}} \leq 1, \quad \textrm{for all $\sigma$},
\]
implies that for all $k, l \in K$,
\begin{equation}
\label{lowbnd2}
\big|
\big\langle \theta(\mu)v,\theta(k)v \big\rangle_{\cV}
-
\big\langle \theta(\mu)v,\theta(l)v \big\rangle_{\cV}
\big|
\leq \big(2\delta + |k^{-1}l|_{S_\delta}\big)\|v\|^2_{\cV},
\end{equation}
where $|\cdot|_{S_\delta}$ is defined as in \eqref{def_normsigmaS}, by splitting the sum in \eqref{lowbnd} into two sums, one over $\sigma$ in $S_\delta^c$ and one over $\sigma$ in $S_\delta$. \\

We note that
\begin{equation}\label{ReIntegral}
  \int_K \real \big\langle \theta(\mu)v,\theta(k)v \big\rangle_{\cV} \, d\mu(k) = \|\theta(\mu)v\|^2_{\cV} \geq \|P_\theta \theta(\mu) v\|^2_{\cV}
= \|P_\theta v\|_{\cV}^2.
\end{equation}
Since the integrand in \eqref{ReIntegral} is a real-valued continuous function on $K$, there is a $k_o \in \supp \mu$ such that
\[
\real \big\langle \theta(\mu)v,\theta(k_o)v \big\rangle_{\cV} \geq \|P_\theta v\|_{\cV}^2.
\]
Hence, by \eqref{lowbnd2},
\begin{eqnarray*}
\real \big\langle \theta(\mu)v,\theta(k_o k)v \big\rangle_{\cV}
&=&
\real \big\langle \theta(\mu)v,\theta(k_o)v \big\rangle_{\cV} - \big(\real \big\langle \theta(\mu)v,\theta(k_o)v \big\rangle_{\cV} - \real \big\langle \theta(\mu)v,\theta(k_o k)v \big\rangle_{\cV} \big) \\
&\geq &
\|P_\theta v\|_{\cV}^2 - \big| \big\langle \theta(\mu)v,\theta(k_o)v \big\rangle_{\cV} - \big\langle \theta(\mu)v,\theta(k_o k)v \big\rangle_{\cV} \big| \\
&\geq &
\|P_\theta v\|_{\cV}^2 - \big(2\delta + |k|_{S_\delta}\big)\|v\|^2_{\cV},
\end{eqnarray*}
for all $k \in K$, which finishes the proof of the first part of Proposition \ref{prop_main2}. \\

For the second part, let us assume that $\mu$ is absolutely continuous with respect to the Haar probability measure $m_K$
on $K$ with Radon-Nikodym derivative $\rho$ and that $\|\rho\|_{L^2(K)} \leq \xi$ for some $\xi$. Since the $L^1$-norm of
$\rho$ is one, we see that $\xi$ must be at least equal to one (with equality if and only if $\rho$ is essentially constant).
By Theorem 15.2.4 in \cite{Dix} (or Section 5.3 in \cite{Foll}), which is an extension of Parseval's Theorem to the setting at hand,
\begin{equation}
\label{eq_parseval}
\|\rho\|^2_{L^2(K)} = \sum_{\sigma \in \widehat{K}} d_\sigma \|\sigma(\mu)\|_{\textrm{HS}}^2,
\end{equation}
where $\widehat{K}$ denotes the set of unitary equivalence classes of strongly continuous \emph{irreducible} unitary $K$-representations; for each $\sigma \in \widehat{K}$, we denote by $d_\sigma$ the (finite) dimension of the corresponding
$K$-representation $(\cW,\sigma)$ and $\|\sigma(\mu)\|_{\textrm{HS}}$ denotes the Hilbert-Schmidt norm of the linear operator $\sigma(\mu)$, explicitly defined as
\[
\|\sigma(\mu)\|_{\textrm{HS}}^2 = \trace \sigma(\mu)^* \sigma(\mu) = \sum_{i=1}^{d_\sigma} \|\sigma(\mu)e_i\|_{\cW}^2,
\]
where $(e_i)$ is some (any) orthonormal basis of $\cW$. By choosing $e_1$ to be the unit vector in $\cW$ for which
$\|\sigma(\mu)\|_{\textrm{op}} = \|\sigma(\mu)e_1\|_{\cW}$, and then extending to a basis, we see that we always have $\|\sigma(\mu)\|_{\textrm{HS}} \geq \|\sigma(\mu)\|_{\textrm{op}}$. By \eqref{eq_parseval}, and the fact that $d_{\id} = 1$ and $\id(\mu) = \id$, we have
\[
1 + \delta^2 d_K |S_\delta| \leq \sum_{\sigma \in \widehat{K}} d_\sigma \|\sigma(\mu)\|^2_{\textrm{op}}
\leq \sum_{\sigma \in \widehat{K}} d_\sigma \|\sigma(\mu)\|^2_{\HS} = \|\rho\|^2_{L^2(K)} \leq \xi^2,
\]
where $d_K = \min_{\sigma \neq \id} d_{\sigma}$, and thus
\[
|S_\delta| \leq \frac{\xi^2-1}{d_K \delta^2},
\]
which finishes the proof of the second part of Proposition \ref{prop_main2}.

\bibliographystyle{amsplain}
\bibliography{BjGr_Bohrtriple3}

\end{document}